\documentclass[psamsfonts,11pt]{amsart}
\usepackage{amssymb}
\usepackage[left=1.29in,right=1.29in,top=1in,bottom=1in]{geometry}
\newtheorem{em-deff}{Definition}[section]
\newtheorem{lemma}[em-deff]{Lemma}
\newtheorem{theorem}[em-deff]{Theorem}
\newtheorem{corollary}[em-deff]{Corollary}

\newtheorem{proposition}[em-deff]{Proposition}
\newtheorem{em-fact}[em-deff]{Fact}
\newtheorem{em-example}[em-deff]{Example}
\newtheorem{claim}[em-deff]{Claim}
\newtheorem{problem}[em-deff]{Problem}

\newtheorem{em-remark}[em-deff]{Remark}
\newtheorem{question}[em-deff]{Question}

\newenvironment{example}{\begin{em-example} \em }{ \end{em-example}}
\newenvironment{remark}{\begin{em-remark} \em }{ \end{em-remark}}
\newenvironment{deff}{\begin{em-deff} \em }{ \end{em-deff}}

\newcommand{\N}{\mathbb N}
\newcommand{\Z}{\mathbb Z}

\def\admissible{Stoyanov}
\def\Min{\mathbf{Min}}
\def\Ps{\mathbf{Ps}}
\def\cf{\mathrm{cf}}
\def\wtd{\mathrm{wtd}}
\def\cont{\mathfrak c}

\title{Minimal pseudocompact group topologies on free abelian groups}

\author[D. Dikranjan]{Dikran Dikranjan}
\address[Dikran Dikranjan]{Universit\`a di Udine, Dipartimento di Matematica e Informatica
\\
 via delle Scienze, 206 - 33100 Udine, Italy}
\email{dikran.dikranjan@dimi.uniud.it}

\author[A. Giordano Bruno]{Anna Giordano Bruno}
\address[Anna Giordano Bruno]{Universit\`a di Udine, Dipartimento di Matematica e Informatica
\\ via delle Scienze, 206 - 33100 Udine, Italy}
\email{anna.giordanobruno@dimi.uniud.it}

\author[D. Shakhmatov]{Dmitri Shakhmatov}
\address[Dmitri Shakhmatov]{Graduate School of Science and Engineering,
Division of Mathematics, Physics and Earth Sciences\\
Ehime University, Matsuyama 790-8577, Japan}
\email{dmitri@dpc.ehime-u.ac.jp}

\dedicatory{Dedicated to Robert Lowen on the occasion of his 60th anniversary}


\thanks{To appear in: {\bf Topology and its Applications\/}.}

\thanks{The first author was partially supported by the SRA grants P1-0292-0101 and J1-9643-0101}

\thanks{The third author was partially supported by the Grant-in-Aid for Scientific 
Research no.~19540092 by the Japan Society for the Promotion of Science (JSPS)}
\begin{document}

\begin{abstract} 
A Hausdorff topological group $G$ is minimal if every continuous isomorphism $f:G\to H$ between $G$ and a Hausdorff topological group $H$ is 
open. Significantly strengthening a 1981 result of Stoyanov, we prove the following theorem: For every infinite minimal 
abelian
group $G$ there exists a sequence $\{\sigma_n:n\in\N\}$ of cardinals such that 
$$
w(G)=\sup\{\sigma_n:n\in\N\}\ \ \mbox{ and } \ \ 
\sup\{2^{\sigma_n}:n\in\N\}\leq |G| \leq 2^{w(G)},
$$
where $w(G)$ is the weight of $G$. If $G$ is an infinite minimal abelian group, then either $|G|=2^\sigma$ for some cardinal $\sigma$, or $w(G)=\min\{\sigma:  |G|\le 2^\sigma\}$;
moreover,
the equality $|G|=2^{w(G)}$ holds whenever $\cf(w(G))> \omega$.  

For a cardinal $\kappa$, we denote by $F_\kappa$ the free abelian group with $\kappa$ many generators. If $F_\kappa$ admits a pseudocompact group topology, then $\kappa\ge\cont$, where $\cont$ is the cardinality of the continuum. We show that the existence of a minimal pseudocompact group topology  on $F_\cont$ is equivalent to the Lusin's Hypothesis $2^{\omega_1}=\mathfrak c$.  For $\kappa>\mathfrak c$, we prove that  $F_\kappa$  admits a (zero-dimensional) minimal pseudocompact group topology if and only if $F_\kappa$ has both a minimal group topology and a pseudocompact group topology.
If $\kappa>\cont$, then $F_\kappa$  admits a connected minimal pseudocompact group topology of weight $\sigma$ if and only if $\kappa=2^\sigma$. Finally, we establish that no infinite torsion-free abelian group can be equipped with a locally connected minimal group topology. 
\end{abstract}

\maketitle

{\em Throughout this paper all topological groups are Hausdorff.\/} We denote by $\mathbb Z$, $\mathbb P$ and $\mathbb N$
respectively the set of integers, the set of primes and  the set of natural numbers. Moreover $\mathbb Q$ denotes the  group of rationals and $\mathbb R$ the group of reals.
For $p\in\mathbb P$ the symbol $\mathbb Z_p$ is used for the group of $p$-adic integers. The symbol $\mathfrak c$ stands for the cardinality of the continuum. For a topological group $G$ the symbol $w(G)$ stands for the weight of $G$.
The Pontryagin dual of a topological abelian group $G$ is denoted by $\widehat G$.  If $H$ is a group and $\sigma$ is a cardinal, then $H^{(\sigma)}$ is used to denote the direct sum of $\sigma$ many copies of the group $H$.
If $G$ and $H$ are groups, then a map $f:G\to H$ is called a {\em monomorphism\/} provided that $f$ is both a group homomorphism and an injection. For undefined terms see \cite{E,Fuchs}.

\begin{deff}
For a cardinal $\kappa$ we use $F_\kappa$ to denote the free abelian group with $\kappa$ many generators.
\end{deff}

\section{Introduction}

The following notion was introduced independently by Choquet (see Do\" \i tchinov \cite{Do}) and Stephenson \cite{S}.

\begin{deff}
A Hausdorff group topology $\tau$ on  a group $G$ is called  \emph{minimal\/}  provided that every Hausdorff group topology $\tau'$ on $G$ such that $\tau'\subseteq \tau$ satisfies $\tau' = \tau$. Equivalently, a Hausdorff topological group $G$ is minimal if  every continuous isomorphism $f:G\to H$ between $G$ and  a Hausdorff topological group $H$ is a topological isomorphism.
\end{deff}

There exist abelian groups which admit no minimal group topologies at all, e.g., the group of rational numbers ${\mathbb Q}$ \cite{P5} or Pr\"ufer's group
${\Z}(p^\infty)$  \cite{DP2}. This suggests the general problem to determine the algebraic structure of the minimal abelian groups, or equivalently, the following

\begin{problem}\label{PM}
{\rm \cite[Problem 4.1]{D}} Describe the abelian groups that admit minimal group topologies.
\end{problem}

Prodanov solved Problem \ref{PM} first for all free abelian groups of finite rank \cite{P}, and later on he improved this result extending it to all cardinals $\leq\cont$ \cite{P5}:

\begin{theorem}\label{Prodanov}
{\rm \cite{P, P5}}
For every cardinal $\kappa\leq\cont$, the group $F_\kappa$ admits minimal group topologies.
\end{theorem}

Since $|F_\kappa|=\omega\cdot \kappa$ for each cardinal $\kappa$, uncountable free abelian groups are determined up to isomorphism by their cardinality. This 
suggests the problem of characterizing the cardinality of minimal abelian groups. The following set-theoretic definition is ultimately relevant to this problem. 

\begin{deff}
\label{definition:of:Stoyanov:cardinals}
\begin{itemize}
\item[(i)]
For infinite cardinals $\kappa$ and $\sigma$  the
 symbol $\Min(\kappa,\sigma)$ denotes the following statement:  There exists a sequence of cardinals $\{\sigma_n: n\in \N\}$ such that 
\begin{equation}\label{min-def}
\sigma=\sup_{n\in\N}\sigma_n\ \  \text{and}\ \  \sup_{n\in\N} 2^{\sigma_n}\leq \kappa \leq 2^{\sigma}.
\end{equation}
We say that the sequence $\{\sigma_n: n\in \N\}$ as above {\em witnesses\/} $\Min(\kappa,\sigma)$.
\item[(ii)]
An infinite cardinal number $\kappa$ satisfying $\Min(\kappa,\sigma)$ for some infinite cardinal $\sigma$ will be called a \emph{\admissible\/} cardinal. 
\item[(iii)] For the sake of convenience,  we add to the class of Stoyanov cardinals also all finite cardinals. 
\end{itemize}
\end{deff}

The cardinals from item (ii) in the above definition were first introduced by Stoyanov in \cite{S3} under the name ``permissible cardinals''.  Their importance is evident from the following fundamental result of Stoyanov providing a complete characterization of the possible cardinalities of minimal abelian groups, thereby solving Problem \ref{PM} for all free abelian groups:

\begin{theorem}\label{Lucho}\label{Lucho2}\emph{\cite{S3}}
\begin{itemize}
	\item[(a)] If $G$ is a minimal abelian group, then $|G|$ is a \admissible\  cardinal.
	\item[(b)] For a cardinal $\kappa$, $F_\kappa$ admits minimal group topologies if and only if $\kappa$ is a \admissible\  cardinal.
\end{itemize}
\end{theorem}

If $\kappa$ is a finite cardinal satisfying (\ref{min-def}), then  $\kappa=2^n$ for some $n\in\N$. On the other hand, every 
finite group is compact and thus minimal. Furthermore, the group $F_n$ admits minimal group topologies for every $n\in\N$ by Theorem \ref{Prodanov}.
It is in order to include also the case of  finite groups in Theorem \ref{Lucho}(a) and finitely generated groups in Theorem \ref{Lucho}(b) that we decided to add item (iii) to  Definition \ref{definition:of:Stoyanov:cardinals}.

It is worth noting that the commutativity of the group in Theorem \ref{Lucho}(b) is important because all restrictions on the cardinality disappear in the case of (non-abelian) free groups:

\begin{theorem}
\emph{\cite{Sh, Remus}}
\label{minimality:of:free:groups}
Every free group admits a minimal group topology.
\end{theorem}

For free groups with infinitely many generators this theorem has been proved in \cite{Sh}. The remaining case was covered in \cite{Remus}.

A subgroup $H$ of a topological group $G$ is \emph{essential\/} 
(in $G$) if  $H\cap N\not=\{e\}$ for every closed normal subgroup $N$ of $G$ with
$N\not=\{e\}$, where $e$ is the identity element of $G$ \cite{P,S}. This notion is a crucial ingredient of the so-called
``minimality criterion'', due to Prodanov and Stephenson \cite{P,S}, describing the dense minimal subgroups of compact groups. 
\begin{theorem} {\rm (\cite{P,S}; see also \cite{DP1,DPS})}
\label{minimality-criterion}
A dense subgroup $H$ of a compact group  $G$ is minimal if and only if $H$ is essential in $G$.
\end{theorem}

A topological group $G$ is \emph{pseudocompact} if every continuous real-valued function defined on $G$ is bounded \cite{H}. In the  spirit of  Theorem \ref{Lucho2}(b)
characterizing the free abelian groups admitting \emph{minimal} topologies, one can also describe the free abelian groups that admit \emph{pseudocompact} 
group topologies (\cite{CReFree,DS}; see Theorem \ref{Ps-free}). The aim of this article is to provide {\em simultaneous\/} minimal and pseudocompact topologization of free abelian groups. To achieve this goal, we need an alternative description of Stoyanov cardinals obtained in Proposition \ref{uniq} as well as a more precise form of Theorem  \ref{Lucho}(a) given in Theorem \ref{Step0}.

We finish this section with a fundamental restriction on the size of pseudocompact groups due to van Douwen.  

\begin{theorem}\label{van:Douwen}
\emph{\cite{vD1}}
If $G$ is an infinite pseudocompact group, then $|G|\ge\cont$.
\end{theorem}

\section{Main results}
\label{results:section}
\subsection{Cardinality and weight of minimal abelian groups}

Let $\kappa$ be a cardinal. Recall that the {\em cofinality\/} $\cf(\kappa)$ of $\kappa$ is defined to be the smallest cardinal $\varkappa$ such that 
there exists a transfinite sequence $\{\tau_\alpha:\alpha\in \varkappa\}$ of cardinals such that $\kappa=\sup\{\tau_\alpha:\alpha\in\varkappa\}$ and 
$\tau_\alpha<\kappa$ for all $\alpha\in \varkappa$. We say that $\kappa$ is \emph{exponential} if $\kappa=2^\sigma$ for some cardinal $\sigma$, and we  call $\kappa$  \emph{non-exponential\/} otherwise. Recall that $\kappa$ is called a {\em strong limit\/} provided that $2^\mu<\kappa$ for every cardinal $\mu<\kappa$. When $\kappa$ is infinite, we define $\log\kappa=\min\{\sigma:\kappa\le 2^\sigma\}$.

We start this section  with a much sharper version of Theorem \ref{Lucho}(a) showing that the weight $w(G)$ of a minimal 
abelian
group $G$ can be taken as   the cardinal $\sigma$ from Definition \ref{definition:of:Stoyanov:cardinals}(ii) witnessing that $|G|$ is a Stoyanov cardinal:

\begin{theorem}\label{Step0}
If $G$ is an infinite minimal abelian group, then $\Min(|G|,w(G))$ holds.
\end{theorem}

This theorem, along with the complete ``internal'' characterization of the Stoyanov cardinals obtained in Proposition \ref{uniq}
permits us to  establish some new important relations between the cardinality and the weight of an arbitrary minimal abelian group. 

\begin{theorem}
\label{weight:minimal:group}
If $\kappa$ is a  cardinal with $\cf(\kappa)>\omega$ and $G$ is a minimal abelian group such that $w(G) \geq \kappa$, then $|G|\geq 2^{\kappa}$. 
\end{theorem}

Let us recall that $|G|=2^{w(G)}$ holds for every compact group $G$ \cite{ComfortTG}. Taking $\kappa=w(G)$ in Theorem \ref{weight:minimal:group} we obtain
the following extension of this property to all minimal abelian groups:

\begin{corollary}
\label{weight:minimal:group:corollary} Let $G$ be a minimal abelian group with $\cf(w(G))>\omega$. Then $|G|=2^{w(G)}$. 
\end{corollary}

Example \ref{divisible:example}(a) below and Theorem \ref{minimality:of:free:groups} show that neither $\cf(w(G))>\omega$ nor ``abelian'' can be removed in Corollary \ref{weight:minimal:group:corollary}. 

Taking  $\kappa= \omega_1$ in Theorem \ref{weight:minimal:group} one obtains the following surprising metrizability criterion for ``small'' minimal abelian groups:

\begin{corollary}\label{metrization}
A minimal abelian group of size $<2^{\omega_1}$ is metrizable. 
\end{corollary}

The condition $\cf(w(G))>\omega$ plays a prominent role in the above results. In particular,  Corollary \ref{weight:minimal:group:corollary} implies that $\cf(w(G))=\omega$
for a minimal abelian group with $|G|<2^{w(G)}$. Our next theorem gives a more precise information in this direction.

\begin{theorem}\label{w=log||}
Let $G$ be an infinite minimal abelian group such that $|G|$ is a non-exponential cardinal. Then $w(G)=\log|G|$ and $\cf(w(G))=\omega$.
\end{theorem}

Under the assumption of GCH, the equality $w(G)=\log|G|$ holds true for every compact group. Theorem 
\ref{w=log||} establishes this property in ZFC for all minimal abelian groups
of non-exponential size.  Let us note that the restraint ``non-exponential'' cannot be omitted, even in the compact case. Indeed,  the equality $w(G)=\log|G|$ may fail for compact abelian groups:  Under the Lusin's Hypothesis $2^{\omega_1}=\cont$,  for the group $G=\Z(2)^{\omega_1}$ one has $w(G)=\omega_1\ne\omega=\log\cont=\log|G|$.

\begin{example}\label{example1}
There exists a consistent example of a compact abelian group $G$ such that $\cf(w(G))=\omega$ and $w(G)>\log |G|$ (see Example  \ref{Easton} (b)). 
\end{example}

\subsection{Minimal pseudocompact group topologies on free abelian groups}

Since pseudocompact metric spaces are compact, from Corollary \ref{metrization} we immediately get the following: 

\begin{corollary}
\label{min:psc:imply:compact:for:small:groups}
Let $G$ be an abelian group such that $|G|<2^{\omega_1}$. Then $G$ admits a minimal pseudocompact group topology if and only if $G$ admits a compact metric group topology. 
\end{corollary}

By Theorem \ref{van:Douwen}, this corollary is vacuously true under the Lusin's Hypothesis $2^{\omega_1}=\mathfrak c$. 

Corollary \ref{min:psc:imply:compact:for:small:groups} shows that for abelian groups of ``small size'' minimal and pseudocompact topologizations are connected in some sense by compactness. We shall see in Corollary \ref{div} below that the same phenomenon happens for divisible abelian groups, irrespectively of their size.

Rather surprisingly, the mere existence of a minimal group topology on $F_\kappa$ quite often implies the existence of a group topology on $F_\kappa$ that is both minimal and pseudocompact. In other words, one often gets pseudocompactness ``for free''.

\begin{theorem}\label{psc:topologies:from:minimal:ones}
Let $\kappa$ and $\sigma$ be   infinite  cardinals. Assume also that $\sigma$ is not a strong limit. If $F_\kappa$ admits a minimal group topology of weight $\sigma$,
then  
$F_\kappa$ also admits a zero-dimensional minimal pseudocompact group topology of weight $\sigma$. 
\end{theorem}

Recall that the {\em beth cardinals\/} $\beth_\alpha$ are defined by recursion on $\alpha$ as follows. Let $\beth_0=\omega$. If $\alpha=\beta+1$ is a successor ordinal,
then $\beth_\alpha=2^{\beth_\beta}$. If $\alpha$ is a limit ordinal, then $\beth_\alpha=\sup\{\beth_\beta:\beta\in\alpha\}$. 

The restriction on weight in 
Theorem \ref{psc:topologies:from:minimal:ones}
is necessary, as our next example demonstrates.

\begin{example}
\label{example:of:minimal:non-pseudocompact}
Let 
$\kappa=\beth_\omega$.
Clearly, the sequence $\{\beth_n:n\in\N\}$ witnesses that $\kappa$ is a Stoyanov cardinal, so {\em $F_\kappa$ admits a minimal group topology\/} $\tau$ by Theorem \ref{Lucho}(b). On the other hand, since $\kappa$ is a strong limit cardinal with $\cf(\kappa)=\omega$ and 
$|F_\kappa|=\kappa$, the group {\em $F_\kappa$ does not admit any pseudocompact group topology\/} by the result of van Douwen \cite{vD1}.
Note that 
$w(F_\kappa,\tau)=\log|F_\kappa|=\log\kappa=\kappa$ 
by 
Theorem \ref{w=log||}, so 
$\sigma=w(F_\kappa,\tau)$ 
is a strong limit cardinal.
\end{example}

``Going in the opposite direction'', in Example \ref{example:of:pseudocompact:non-minimal} below we will define a cardinal $\kappa$ such that $F_\kappa$ admits a pseudocompact group topology of weight $\sigma$ that is not a strong limit cardinal, and yet $F_\kappa$ does not admit any minimal group topology. These two examples show that the existence of a minimal group topology and the existence of a pseudocompact group topology on a free abelian group are ``independent events''. 

For a free group of size $>\cont$ that admits both a minimal group topology and a pseudocompact group topology, the next theorem discovers the surprising possibility of ``simultaneous topologization'' with a topology which is both minimal and pseudocompact. Moreover, it turns out that  this topology can  also be chosen to be zero-dimensional. 

\begin{theorem}\label{MPs}
For every cardinal $\kappa>\mathfrak c$ the following conditions are equivalent:
  \begin{itemize}
    \item[(a)] $F_\kappa$ admits both a minimal group topology and a pseudocompact group topology;
    \item[(b)] $F_\kappa$ admits a minimal pseudocompact group topology;
    \item[(c)] $F_\kappa$ admits a zero-dimensional minimal pseudocompact group topology.
	\end{itemize}
\end{theorem}

The free abelian group group $F_\mathfrak c$ of cardinality $\mathfrak c$ admits a minimal group topology (Theorem \ref{Prodanov}) and a pseudocompact group topology \cite{DS}.  Our next theorem shows that the statement ``$F_\mathfrak c$ admits a minimal pseudocompact group topology'' is both consistent with and independent of ZFC.

\begin{theorem}\label{continuum}
The following conditions are equivalent:
\begin{itemize}
	\item[(a)] $F_\cont$ admits a minimal pseudocompact group topology;
	\item[(b)] $F_\mathfrak c$ admits a connected minimal pseudocompact group topology;
	\item[(c)] $F_\mathfrak c$ admits a zero-dimensional minimal pseudocompact group topology;
	\item[(d)] the Lusin's Hypothesis $2^{\omega_1}=\mathfrak c$ holds.
\end{itemize}
\end{theorem}

Since every infinite pseudocompact group has cardinality $\geq \mathfrak c$ (Theorem \ref{van:Douwen}), 
Theorems \ref{MPs} and \ref{continuum} provide a complete description of free abelian groups that have a minimal (zero-dimensional) pseudocompact group topology. The equivalence of (a) and (b) in Theorem \ref{MPs} (respectively, (a) and (d) in Theorem \ref{continuum}) was announced without proof in \cite[Theorem 4.11]{D}.

Motivated by Theorem \ref{MPs}(c)  and Theorem \ref{continuum}(c), where the minimal pseudocompact topology can be additionally 
chosen
zero-dimensional (or connected,
in Theorem \ref{continuum}(b)), we arrive at the following natural question: {\em  If $\kappa$ is a cardinal such that $F_\kappa$ admits a minimal group topology $\tau_1$ and  a pseudocompact group topology $\tau_2$, and one of these topologies is connected, does then $F_\kappa$ admit a connected minimal pseudocompact group topology $\tau_3$?\/} Theorem \ref{continuum} answers  this question  in the case of $F_\cont$. The next theorem gives an answer for $\kappa>\mathfrak c$, showing a symmetric behavior, as far as connectedness is concerned. This should be compared with the  equivalent items in Theorem \ref{continuum}
where item (a) contains no restriction beyond minimality and pseudocompactness, whereas item (c) contains ``zero-dimensional''.

\begin{theorem}\label{connminpsc}
Let $\kappa$ and $\sigma$ be infinite cardinals with $\kappa>\mathfrak c$. The following conditions are equivalent:
\begin{itemize}
	\item[(a)]$F_\kappa$ admits a connected minimal pseudocompact group topology (of weight $\sigma$);
	\item[(b)]$F_\kappa$ admits a connected minimal group topology (of weight $\sigma$);
	\item[(c)]$\kappa$ is exponential ($\kappa=2^\sigma$).
\end{itemize}
\end{theorem}

This theorem is ``asymmetric'' in some sense toward minimality. Indeed, item (b) should be compared with the fact that the existence of a connected pseudocompact group topology on $F_\kappa$ need not necessarily imply that $F_\kappa$ admits a connected minimal group topology (see Example \ref{Example:LaST}).

If a free abelian group admits a pseudocompact group topology, then it admits also a pseudocompact group topology which is both connected and locally connected \cite[Theorem 5.10]{DS}. When minimality is added to the mix, the situation becomes totally different. In Example \ref{Example:LaST} below we exhibit a free abelian group $F_\kappa$ that admits a connected, locally connected, pseudocompact group topology, and yet  
$F_\kappa$ does not have any connected minimal group topology.
Even more striking is the following
\begin{theorem}
\label{no:locally:connected}
A locally connected minimal torsion-free abelian group is trivial.
\end{theorem}
Theorem \ref{no:locally:connected} strengthens significantly \cite[Corollary 8.8]{DS} by replacing ``compact'' in it with ``minimal''.
\begin{corollary}
No free abelian group admits a locally connected, minimal 
group topology.
\end{corollary}

The reader may wish to compare this corollary with Theorems \ref{continuum} and \ref{connminpsc}.

The paper is organized as follows.  In Section \ref{admissible:section} we give some properties of Stoyanov cardinals, while Section \ref{section:4} contains all necessary facts concerning pseudocompact topologization. The culmination here is Corollary \ref{Simultaneous} establishing that, roughly speaking, if  $F_\kappa$ admits a minimal group topology $\tau_1$ and a pseudocompact group topology $\tau_2$, then one can assume, without loss of generality, that this pair satisfies  $w(F_\kappa, \tau_1)=w(F_\kappa, \tau_2)$.   Sections \ref{section:5} and \ref{section:6} prepare the remaining necessary tools for the proof of the  main results, deferred  to Section \ref{proofs:section}. 
Finally, in Section \ref{section:8} we discuss the counterpart of the simultaneous minimal and pseudocompact topologization 
for other classes of abelian groups such as divisible groups, torsion-free groups and torsion groups, as well as the same problem for (non-commutative) free groups. 

\section{Properties of \admissible\ cardinals}
\label{admissible:section}

We start with an example of small  \admissible\ cardinals. 

\begin{example}\label{example:small:Stoyanov}
If $\omega\leq\kappa\leq\mathfrak c$, then $\Min(\kappa,\omega)$. \end{example}
In our next example we discuss the connection between $\Min(\kappa,\sigma)$ and the property of $\kappa$ to be exponential. 

\begin{example}\label{exp-non-exp}
Let $\kappa$ be an infinite cardinal.
\begin{itemize}
	\item[(a)] 
{\em If $\kappa=2^{\sigma}$,
then $\Min(\kappa,\sigma)$ holds. In particular,
an exponential cardinal is \admissible.\/}
	\item[(b)] {\em  If $\{\sigma_n:n\in\N\}$ is a sequence of cardinals witnessing $\Min(\kappa,\sigma)$
such that $\sigma=\sigma_m$ for some $m\in\N$, then $\kappa=2^\sigma$\/}.  Indeed,  (\ref{min-def}) and our assumption yield
$$
2^\sigma=2^{\sigma_m}\le \sup_{n\in\N}2^{\sigma_n}\le \kappa\le 2^\sigma.
$$
Hence $\kappa=2^\sigma$.

	\item[(c)] \emph{If $\cf(\sigma)>\omega$, then $\Min(\kappa,\sigma)$ if and only if $\kappa=2^\sigma$}. If $\kappa=2^\sigma$, then $\Min(\kappa,\sigma)$ by item (a). Assume $\Min(\kappa,\sigma)$, and let $\{\sigma_n:n\in\N\}$ be a sequence of cardinals 
witnessing $\Min(\kappa,\sigma)$.
From (\ref{min-def}) and $\cf(\sigma)>\omega$ 
we get $\sigma=\sigma_m$ for some $m\in\N$. Applying  item (b) gives $\kappa=2^\sigma$.  
\end{itemize}
\end{example}

Clearly, $\Min(\kappa,\sigma)$ implies  $\sigma\geq \log\kappa$. We show now that this inequality becomes an equality in case $\kappa$ is non-exponential. 

\begin{lemma}\label{non-exp}
Let $\kappa$ be a non-exponential infinite cardinal. Then:
\begin{itemize}
	\item[(a)]$\Min(\kappa,\sigma)$ if and only if $\cf(\sigma)=\omega$ and $\log\kappa=\sigma$;
	\item[(b)] $\Min(\kappa,\log\kappa)$ if and only if $\cf(\log\kappa)=\omega$.
\end{itemize}
\end{lemma}
\begin{proof}
(a) To prove the ``only if'' part, assume that $\Min(\kappa,\sigma)$ holds,  and let $\{\sigma_n:n\in\N\}$ be a sequence of cardinals witnessing
$\Min(\kappa,\sigma)$.
Since $\kappa\le 2^\sigma$ by (\ref{min-def}), we have 
$\log\kappa\le \sigma$. Assume $\log \kappa<\sigma$.  From  (\ref{min-def}) we conclude that
$\log \kappa\le\sigma_m$ for some $m\in\N$. Therefore 
$$
2^{\log \kappa}\leq2^{\sigma_m}\leq \sup_{n\in\N}2^{\sigma_n}\leq\kappa\le 2^{\log \kappa}
$$ 
by (\ref{min-def}). Thus $\kappa=2^{\log \kappa}$ is an exponential cardinal, a contradiction. This proves that $\sigma=\log \kappa$.

To prove the ``if'' part, assume that $\cf(\sigma)=\omega$ and $\log\kappa=\sigma$. Then  there exists a sequence of cardinals $\{\sigma_n:n\in\N\}$ such that $\sigma=\sup_{n\in\N}\sigma_n$ and $\sigma_n<\sigma=\log\kappa$ for every $n\in\N$. In particular,
$2^{\sigma_n}<\kappa$ for every $n\in\N$. Consequently, 
$$
\sup_{n\in\N}2^{\sigma_n}\leq\kappa\leq 2^{\log\kappa}=2^{\sigma}.
$$
That is,  (\ref{min-def}) holds. Therefore, the sequence $\{\sigma_n:n\in\N\}$ witnesses
$\Min(\kappa,\sigma)$.

\medskip

Item (b) follows from item (a). 
\end{proof}

\begin{example}\label{Easton}
Let $\kappa$ and $\sigma$ be cardinals. According to Example \ref{exp-non-exp}(a),  $\Min(\kappa,\sigma)$ does not imply $\cf(\sigma)=\omega$ in case $\kappa$ is exponential. (Indeed, it suffices to take $\kappa=2^\sigma$ with $\cf(\sigma)>\omega$.) 
\begin{itemize}
    \item[(a)]  Let us show that the assumption
 ``$\kappa$ is non-exponential''  in Lemma \ref{non-exp}(a) is necessary (to prove that $\Min(\kappa,\sigma)$
    implies $\log\kappa=\sigma$) even in the case $\cf(\sigma)=\omega$. To this end, use an appropriate Easton model \cite{Easton} satisfying  
$$
2^{\omega_{\omega+1}}=\omega_{\omega+2}\ \ \mbox{ and }\ \  2^{\omega_n}=\omega_{\omega+2}\mbox{ for all }n\in\N.
$$ 
  Let $\kappa=\omega_{\omega+2}$ and $\sigma=\omega_\omega$. Then $2^\sigma=\kappa$ as $2^{\omega_{\omega+1}}=2^{\omega_n}=\kappa$ for every $n\in\N$.
  So $\Min(\kappa,\sigma)$ holds by Example \ref{exp-non-exp}(a). Moreover $\cf(\sigma)=\omega$ and $\log\kappa=\omega_0<\omega_\omega=\sigma$. 
\item[(b)] Using the cardinals $\kappa$ and $\sigma$ from item (a) we can give now the example anticipated in Example \ref{example1}. 
Let $G=\Z(2)^\sigma$. Then $w(G) = \sigma$, so $\cf(w(G))=\omega$ and  yet $\log |G|= \log 2^\sigma=\log \kappa=\omega_0<\sigma=w(G)$.
\end{itemize}
\end{example}

The next proposition, summarizing the above results, provides an alternative description of the infinite Stoyanov cardinals that makes no use of the somewhat ``external'' condition (\ref{min-def}).

\begin{proposition}\label{uniq} Let  $\kappa$  be an infinite cardinal. 
\begin{itemize}
	\item[(a)] If $\kappa$  is exponential, then $\Min(\kappa,\sigma)$ holds for every cardinal $\sigma$ with $\kappa=2^{\sigma}$.
	\item[(b)] If $\kappa$ is non-exponential, then $\Min(\kappa,\sigma)$ is equivalent to $\sigma=\log\kappa$ and $\cf(\log\kappa)=\omega$.
\end{itemize}
\end{proposition}

\begin{proof} Item (a) follows from Example \ref{exp-non-exp}(a), and item  (b) follows from Lemma \ref{non-exp}(a).  
\end{proof}

\section{Cardinal invariants related to pseudocompact groups}
\label{section:4}

Recall that a subset $Y$ of a space $X$ is said to be {\em $G_\delta$-dense in $X$\/} provided that $Y\cap B\not=\emptyset$ for every non-empty $G_\delta$-subset $B$ of $X$.

The following theorem describes pseudocompact groups in terms of their completion.

\begin{theorem}\label{cr-theorem}\emph{\cite[Theorem 4.1]{CR}}
A precompact group $G$ is pseudocompact if and only if  $G$ is $G_\delta$-dense in its 
completion.
\end{theorem}

\begin{deff}
\label{Ps-def}
\begin{itemize}
\item[(i)]
If $X$ is a non-empty set and $\sigma$ is an infinite cardinal, then a set $F\subseteq X^\sigma$ is \emph{$\omega$-dense} in $X^\sigma$, provided that for every countable set $A\subseteq \sigma$ and each function $\varphi\in X^A$ there exists $f\in F$ such that $f(\alpha)=\varphi(\alpha)$ for all $\alpha\in A$.
\item[(ii)]
If $\kappa$ and $\sigma\geq\omega$ are cardinals, then $\Ps(\kappa,\sigma)$ abbreviates the sentence ``there exists an $\omega$-dense set $F\subseteq\{0,1\}^\sigma$ with $|F|=\kappa$''.
\item[(iii)]
For an infinite cardinal $\sigma$ let $m(\sigma)$ denote the minimal cardinal $\kappa$ such that $\Ps(\kappa,\sigma)$ holds.
\end{itemize}
\end{deff}

Items (i) and (ii) of the above definition are taken from \cite{CEG} except for the notation $\Ps(\kappa,\sigma)$ that appears in \cite[Definition 2.6]{DS}.
Item (iii) is equivalent to the definition of the cardinal function $m(\sigma)$ of Comfort and Robertson \cite{CRob1}. It is worth noting that 
$m(\sigma)=\delta(\sigma)$ for every infinite cardinal $\sigma$, where $\delta(-)$ is the cardinal function defined by Cater, Erd\"os and Galvin \cite{CEG}.

The set-theoretical  condition $\Ps(\kappa,\sigma)$ is ultimately related to the existence of pseudocompact group topologies. 

\begin{theorem}\label{CRThm}
\emph{(\cite{CRob1}; see also \cite[Fact 2.12 and Theorem 3.3(i)]{DS})}\label{Ps-top}
Let $\kappa$ and $\sigma\geq\omega$ be cardinals. Then $\Ps(\kappa,\sigma)$ holds if and only if there exists a group $G$ of cardinality $\kappa$ which admits a pseudocompact group topology of weight $\sigma$.
\end{theorem}

Moreover, the condition $\Ps(\kappa,\sigma)$ completely describes free abelian groups that admit pseudocompact group topologies.

\begin{theorem}\emph{(\cite{CReFree}, \cite[Theorem 5.10]{DS})}\label{Ps-free}
If $\kappa$ is a cardinal, then $F_\kappa$ admits a pseudocompact group topology of weight $\sigma$
if and only if $\Ps(\kappa,\sigma)$ holds.
\end{theorem}

In the next lemma we summarize some properties of the cardinal function $m(-)$  for future reference.

\begin{lemma}\label{m-log}\emph{(\cite{CEG}; see also \cite[Theorem 2.7]{CRob1})} Let $\sigma$ be an infinite cardinal. Then:
\begin{itemize}
	\item[(a)] $m(\sigma)\geq 2^\omega$ and $\cf(m(\sigma))>\omega$;
	\item[(b)] $\log\sigma\leq m(\sigma)\leq(\log \sigma)^\omega$;
	\item[(c)] $m(\lambda)\leq m(\sigma)$ whenever $\lambda$ is a cardinal with $\lambda\leq\sigma$.
\end{itemize}
\end{lemma}

Some useful properties of the condition $\Ps(\lambda,\kappa)$ are collected in the next proposition. Items (a) and (b) are part of \cite[Lemmas 2.7 and 2.8]{DS}, and items (d) and (e) are particular cases of \cite[Lemma 3.4(i)]{DS}.

\begin{proposition}\label{DS}
\begin{itemize}
	\item[(a)] $\Ps(\mathfrak c,\omega)$ holds, and moreover, $m(\omega)=\mathfrak c$; also $\Ps(\mathfrak c,\omega_1)$ holds.
	\item[(b)] If $\Ps(\kappa,\sigma)$ holds for some cardinals $\kappa$ and
$\sigma\geq\omega$, then $\kappa\geq\mathfrak c$, and $\Ps(\kappa',\sigma)$ holds for every cardinal $\kappa'$ such that $\kappa\leq\kappa'\leq 2^\sigma$.
	\item[(c)] For cardinals $\kappa$ and $\sigma\geq\omega$, $\Ps(\kappa,\sigma)$ holds if and only if $m(\sigma)\leq \kappa\leq 2^\sigma$.
	\item[(d)] $\Ps(2^\sigma,\sigma)$ and $\Ps\left(2^\sigma,2^{2^\sigma}\right)$ hold for every infinite cardinal $\sigma$.
	\item[(e)] If $\sigma$ is a cardinal such that $\sigma^\omega=\sigma$, then $\Ps(\sigma, 2^\sigma)$ holds.
\end{itemize}
\end{proposition}

\begin{example}
\label{example:of:pseudocompact:non-minimal}
Let $\kappa=\beth_{\omega_1}$ (see the text preceding Example \ref{example:of:minimal:non-pseudocompact} for the definition of $\beth_{\omega_1}$).
One can easily see that 
$\kappa$ 
is not a Stoyanov cardinal (this was first noted by Stoyanov himself). Therefore, the group {\em $F_\kappa$ does not admit any minimal group topology\/} by Theorem \ref{Lucho}(a).
On the other hand, $\kappa=\kappa^\omega$ and Proposition \ref{DS}(e) yield that $\Ps(\kappa,2^\kappa)$ holds. Applying Theorem \ref{Ps-free} we conclude that {\em $F_\kappa$ admits a pseudocompact group topology of weight $2^\kappa$\/}. In particular, $\sigma=2^\kappa$ is not a strong limit.
\end{example}

Example \ref{example:of:pseudocompact:non-minimal} should be compared with Theorem \ref{psc:topologies:from:minimal:ones} where we show that  if $F_\kappa$ admits a minimal group topology of weight $\sigma$ and $\sigma$ is not a a strong limit, then $F_\kappa$ admits also a pseudocompact group topology of weight $\sigma$.

\begin{example}\label{Example:LaST}
Let $\kappa$ be a non-exponential cardinal with $\kappa=\kappa^\omega$ (e.g., a strong limit cardinal of uncountable cofinality). Then, according to Proposition \ref{DS}(e), $\Ps(\kappa, 2^\kappa)$ holds. Therefore {\em $F_\kappa$ admits a pseudocompact group topology (of weight $2^\kappa$) that is both connected and locally connected\/} \cite[Theorem 5.10]{DS}. By Theorem \ref{connminpsc}, {\em $F_\kappa$ does not admit a connected minimal group topology\/} as  $\kappa$ is non-exponential. 
\end{example}

\begin{lemma}
\label{lemma:4:8}
If $\kappa$ and $\sigma$ are infinite cardinals such that $\sigma$ is not a strong limit cardinal,  then $\Min(\kappa,\sigma)$ implies $\Ps(\kappa,\sigma)$.
\end{lemma}

\begin{proof}
Assume that $\Min(\kappa,\sigma)$ holds, and let $\{\sigma_n:n\in\N\}$ be a sequence of cardinals witnessing $\Min(\kappa,\sigma)$.
Since $\sigma$ is not a strong limit cardinal, there exists a cardinal $\mu<\sigma$ such that $\sigma\le 2^\mu$. Since $\sigma=\sup_{n\in\N} \sigma_n$ by (\ref{min-def}), $\mu\le \sigma_n$ for some $n\in\N$. Then  $\sigma\le 2^\mu\le 2^{\sigma_n}$, and so $\log\sigma\le \sigma_n$.  Applying Lemma \ref{m-log}(b)  and (\ref{min-def}), we obtain
$$
m(\sigma)\leq (\log\sigma)^\omega\leq\sigma_n^\omega\leq2^{\sigma_n}\leq\kappa\le 2^\sigma.
$$ 
Hence $\Ps(\kappa,\sigma)$ holds by Proposition \ref{DS}(c).
\end{proof}

Our next example demonstrates that the restriction on the cardinal $\sigma$ in Lemma \ref{lemma:4:8} is necessary.

\begin{example}
Let $\kappa$ be the Stoyanov cardinal from Example \ref{example:of:minimal:non-pseudocompact}. From calculations in that example
one concludes that $\Min(\kappa,\kappa)$ holds. As was shown in Example \ref{example:of:minimal:non-pseudocompact}, $F_\kappa$ does not admit any pseudocompact group topology. Therefore,  $\Ps(\kappa,\sigma)$ fails for every cardinal $\sigma$
(Theorem \ref{Ps-free}).
\end{example}

In the next lemma we show that, if $\kappa$ is a Stoyanov cardinal satisfying $\Ps(\kappa,\lambda)$ for some $\lambda$, then $\Ps(\kappa,\sigma)$ holds also for the cardinal $\sigma$ witnessing that $\kappa$ is Stoyanov.

\begin{lemma}\label{Step1}
Let $\kappa$ and $\sigma$ be infinite cardinals satisfying $\Min(\kappa,\sigma)$. If $\Ps(\kappa,\lambda)$ holds for some infinite cardinal $\lambda$, then $\Ps(\kappa,\sigma)$ holds as well.
\end{lemma}
\begin{proof} By Lemma \ref{lemma:4:8}, it suffices only to consider the case when $\sigma$ is a strong limit cardinal. Let $\{\sigma_n:n\in\N\}$ be a sequence of cardinals witnessing $\Min(\kappa,\sigma)$. If $\sigma=\sigma_n$ for some $n\in\N$, then $\kappa=2^{\sigma}$ by Example \ref{exp-non-exp}(b). Since $\Ps(2^{\sigma},\sigma)$ holds by Proposition \ref{DS}(d), we are done in this case. Suppose now that $\sigma>\sigma_n$ for every $n\in\N$. Since  $\Ps(\kappa,\lambda)$ holds, from Proposition \ref{DS}(c) we get $m(\lambda)\leq\kappa\leq 2^\lambda$. If $\lambda<\sigma$, then $2^\lambda<\sigma$ and so $\kappa<\sigma$. From (\ref{min-def}) we get $\kappa<\sigma_n$ for some $n\in\N$, and then $\sigma_n<2^{\sigma_n}\le \kappa$, a contradiction. Hence $\sigma\leq\lambda$. By Lemma \ref{m-log}(c) $m(\sigma)\leq m(\lambda)\leq \kappa$. Moreover $\kappa\leq 2^\sigma$ by  (\ref{min-def}). It now follows from Proposition \ref{DS}(c)  that $\Ps(\kappa,\sigma)$ holds.
\end{proof}

\begin{corollary}\label{Simultaneous} 
Let $\kappa$ be a non-zero cardinal. If $F_\kappa$ admits a minimal group topology $\tau_1$ and a pseudocompact group topology $\tau_2$, then $F_\kappa$ admits also a pseudocompact group topology $\tau_3$ with $w(F_\kappa, \tau_1)=w(F_\kappa, \tau_3)$.
\end{corollary}

\begin{proof} From Theorem \ref{van:Douwen} we get $\kappa\ge\cont$. Define $\sigma=w(F_\kappa, \tau_1)$.  Clearly, $\sigma$ is infinite.
Applying Theorem \ref{Step0}, we conclude that $\Min(\kappa,\sigma)$ holds. Theorem \ref{Ps-free} yields that $\Ps(\kappa,\lambda)$ holds, where $\lambda=w(F_\kappa,\tau_2)$.
Clearly, $\lambda$ is infinite. Then $\Ps(\kappa,\sigma)$ holds by Lemma \ref{Step1}.  Finally, applying Theorem \ref{Ps-free} once again, we obtain that $F_\kappa$ must admit a pseudocompact group topology $\tau_3$ such that  $w(F_\kappa, \tau_3)=\sigma$. 
\end{proof}

The proof of Corollary \ref{Simultaneous} relies on Theorem \ref{Step0}, which is proved later in Section \ref{proofs:section}. Nevertheless, this does not create any problems, because Corollary \ref{Simultaneous}
is never used thereafter.

\section{Building $G_\delta$-dense $\mathcal{V}$-independent subsets in products}
\label{section:5}

A \emph{variety of groups} $\mathcal V$ is a class of abstract groups closed under subgroups, quotients and products. For a variety $\mathcal V$ and $G\in\mathcal V$ a subset $X$ of $G$ is \emph{$\mathcal V$-independent} if the subgroup $\langle X\rangle$ of $G$ generated by $X$ belongs to $\mathcal{V}$ and for each map $f:X\to H\in\mathcal V$ there exists a unique homomorphism $\overline{f}:\langle X\rangle\to H$ extending $f$. Moreover, the \emph{$\mathcal V$-rank of $G$} is 
$$
r_\mathcal V(G):=\sup\{|X|:X\text{ is a $\mathcal V$-independent subset of }G\}.
$$
In particular, if $\mathcal A$ is the variety of all abelian groups,  then the $\mathcal A$-rank is the usual free rank $r(-)$, and for the variety 
$\mathcal A_p$ of all abelian groups of exponent $p$ (for a prime $p$) the $\mathcal A_p$-rank is the usual $p$-rank $r_p(-)$.

\smallskip
Our first lemma is a generalization of \cite[Lemma 4.1]{DS} that is in fact equivalent to \cite[Lemma 4.1]{DS} (as can be seen from its proof below).
\begin{lemma}
\label{big:rank:of:products}
Let $\mathcal{V}$ be a variety of groups and $I$ an infinite set. For every $i\in I$ let $H_i$ be a group such that $r_{\mathcal{V}}(H_i)\ge \omega$.
Then $r_{\mathcal{V}}\left(\prod_{i\in I} H_i\right)\ge 2^{|I|}$.
\end{lemma}

\begin{proof} 
Define $N=\mathbb{N}\setminus\{0\}$. For every $n\in N$, let $F_n$ be the free group in the variety $\mathcal V$ with $n$ generators.
Define $H=\prod_{n\in N} F_n$, and note that $r_{\mathcal{V}}(H)\ge\omega$.
Since $I$ is infinite, there exists a bijection $\xi: I\times N\to I$. For $(i,n)\in I\times N$, fix a subgroup
$F_{in}$ of $H_{\xi(i,n)}$ isomorphic to $F_n$ (this can be done because
$r_{\mathcal{V}}(H_{\xi(i,n)})\ge \omega$). Then
$\prod_{(i,n)\in I\times N} F_{in}$
is a subgroup of the group
$\prod_{(i,n)\in I\times N} H_{\xi(i,n)}\cong \prod_{i\in I} H_i$,
where $\cong$ denotes the isomorphism between groups.
Clearly,
$$
\prod_{(i,n)\in I\times N} F_{in}\cong \prod_{i\in I}\prod_{n\in  N} F_{in}
\cong \prod_{i\in I}\prod_{n\in  N} F_{n}
\cong 
\prod_{i\in I} H\cong H^I,
$$
so there exists a monomorphism 
$f:H^I\to \prod_{i\in I} H_i$.
Now
$$
r_{\mathcal{V}}\left(\prod_{i\in I} H_i\right)\ge r_{\mathcal{V}} \left(f\left(H^I\right)\right)=r_{\mathcal{V}} \left(H^I\right) \ge 2^{|I|},
$$
 where the the first inequality follows from \cite[Corollary 2.5]{DS} and the last inequality has been proved in \cite[Lemma 4.1]{DS}.
\end{proof}

\begin{lemma}
\label{Ps:in:products}
Suppose that $I$ is an infinite set and $H_i$ is a separable metric space for every $i\in I$.  If $\Ps(\kappa, |I|)$ holds, then the product $H=\prod_{i\in I} H_i$ contains a $G_\delta$-dense subset of size at most $\kappa$.
\end{lemma}

\begin{proof} Let $i\in I$. Since $H_i$ is a separable metric space, $|H_i|\le \cont$, and so we can fix a surjection $f_i:\mathbb{R}\to H_i$.

Let $\theta:\mathbb{R}^I\to H$ be the map defined by  $\theta(g)=\{f_i(g(i))\}_{i\in I}\in H$  for every $g\in \mathbb{R}^I$. Since $\Ps(\kappa, |I|)$ holds, \cite[Lemma 2.9]{DS} allows us to conclude that $\mathbb{R}^I$ contains an $\omega$-dense subset $X$ of size $\kappa$. Define $Y=\theta(X)$. Then $|Y|\le |X|=\kappa$. It remains only to show that $Y$ is $G_\delta$-dense in $H$. Indeed, let $E$ be a non-empty $G_\delta$-subset of $H$. Then there exist a countable subset $J$ of $I$ and $h\in \prod_{j\in J} H_j$ such that $\{h\}\times \prod_{i\in I\setminus J} H_i\subseteq E$. For every $j\in J$ select $r_j\in \mathbb R$ such that $f_j(r_j)=h(j)$. Since $X$ is $\omega$-dense in $\mathbb{R}^I$, there exists $x\in X$  such that $x(j)=r_j$ for every $j\in J$. Now 
\begin{align*}
\theta(x)&=\{f_i(x(i))\}_{i\in I}=\{f_j(x(j))\}_{j\in J}\times \{f_i(x(i))\}_{i\in I\setminus J}\\
&=\{h(j)\}_{j\in J}\times \{f_i(x(i))\}_{i\in I\setminus J}\in \{h\}\times \prod_{i\in I\setminus J} H_i\subseteq E.
\end{align*}
Therefore, $\theta(x)\in Y\cap E\not=\emptyset$.
\end{proof}

\begin{lemma}
\label{set:theoretic:lemma:1}
Let $\kappa\ge\omega_1$ be a cardinal and $G$ and $H$ be topological groups in a variety $\mathcal V$ such that:
\begin{itemize}
\item[(a)] $r_\mathcal V(H)\ge\kappa$,
\item[(b)] $H^\omega$ has a $G_\delta$-dense subset of size at most $\kappa$,
\item[(c)] $G$ has a $G_\delta$-dense subset of size at most $\kappa$.
\end{itemize}
Then $G\times H^{\omega_1}$ contains a $G_\delta$-dense $\mathcal V$-independent subset of size  $\kappa$. 
\end{lemma}

\begin{proof}
Since $\kappa\ge \omega_1$, we have $|\kappa\times\omega_1|=\kappa$, and so we can use item (a) to fix a faithfully indexed $\mathcal V$-independent subset $X=\{x_{\alpha\beta}:\alpha\in\kappa,\beta\in\omega_1\}$ of $H$.   For every $\beta\in\omega_1\setminus\omega$ the topological groups $G\times H^\omega$ and $G\times H^\beta$ are isomorphic, so we  can use items (b) and (c) to fix $\{g_{\alpha\beta}:\alpha\in\kappa\}\subseteq G$ and $\{y_{\alpha\beta}:\alpha\in\kappa\}\subseteq H^\beta$ such that $Y_\beta=\{(g_{\alpha\beta}, y_{\alpha\beta}) : \alpha\in\kappa\}$ is a $G_\delta$-dense subset  of $G\times H^\beta$. 

For $\alpha\in\kappa$ and $\beta\in\omega_1\setminus\omega$ define $z_{\alpha\beta}\in H^{\omega_1}$ by
\begin{equation}
\label{define:z}
z_{\alpha\beta}(\gamma)=
\left\{ \begin{array}{ll}
y_{\alpha\beta}(\gamma), & \mbox{ for } \gamma\in\beta  \\
x_{\alpha\beta}, & \mbox{ for } \gamma\in\omega_1\setminus\beta
\end{array} \right.
\hskip25pt \mbox{ for } \gamma\in\omega_1.
\end{equation}
Finally, define  
$$
Z=\{(g_{\alpha\beta}, z_{\alpha\beta}): \alpha\in\kappa, \beta\in\omega_1\setminus\omega\}\subseteq G\times H^{\omega_1}.
$$

\begin{claim}
$Z$ is $G_\delta$-dense in $G\times H^{\omega_1}$.
\end{claim}

\begin{proof}
Let $E$ be a non-empty $G_\delta$-subset of $G\times H^{\omega_1}$. Then there exist $\beta\in\omega_1\setminus\omega$ and a non-empty $G_\delta$-subset $E'$ of $G\times H^\beta$ such that 
\begin{equation}
\label{Gdelta:equation}
E'\times H^{\omega_1\setminus\beta}\subseteq E.
\end{equation}
Since $Y_\beta$ is $G_\delta$-dense in $G\times H^\beta$, there exists $\alpha\in\kappa$ such that $(g_{\alpha\beta},y_{\alpha\beta})\in E'$. From (\ref{define:z}) it follows that
$z_{\alpha\beta}\restriction_\beta=y_{\alpha\beta}$. Combining this with (\ref{Gdelta:equation}), we conclude that
$(g_{\alpha\beta}, z_{\alpha\beta})\in E$. Thus $(g_{\alpha\beta}, z_{\alpha\beta})\in E\cap Z\not=\emptyset$.
\end{proof}

\begin{claim}
$Z$ is $\mathcal V$-independent.
\end{claim}

\begin{proof}
Let $F$ be a non-empty finite subset of  $\kappa\times(\omega_1\setminus\omega)$. Define 
\begin{equation}
\label{equation:gamma}
\gamma=\max\{\beta\in\omega_1\setminus\omega:\exists\ \alpha\in \kappa\ 
(\alpha,\beta)\in F\}.
\end{equation} 
From (\ref{define:z}) and (\ref{equation:gamma}) it follows that 
$z_{\alpha\beta}(\gamma)=x_{\alpha\beta}$
for all $(\alpha,\beta)\in F$. Therefore, 
$$
X_F=\{z_{\alpha\beta}(\gamma):(\alpha,\beta)\in F\} = \{x_{\alpha\beta}:(\alpha,\beta)\in F\} \subseteq X.
$$
Since $X$ is a $\mathcal V$-independent subset of $H$,  so is $X_F$ \cite[Lemma 2.3]{DS}.
Let $f:G\times H^{\omega_1}\to H$ be the projection homomorphism defined by 
$f(g,h)=h(\gamma)$ for $(g,h)\in G\times H^{\omega_1}$.
Define 
$$
S_F=\{(g_{\alpha\beta},z_{\alpha\beta}):(\alpha,\beta)\in F\}.
$$
Since $G\in\mathcal{V}$, $H\in\mathcal{V}$, $\langle S_F\rangle $ is a subgroup of $G\times H^{\omega_1}$ and $\mathcal{V}$ is a variety, 
$\langle S_F\rangle\in\mathcal{V}$. Since 
$f\restriction_{S_F}: S_F\to H$ 
is an injection
and $f(S_F)=X_F$ is a $\mathcal{V}$-independent subset of $H$, from  \cite[Lemma 2.4]{DS} we obtain that $S_F$ is $\mathcal{V}$-independent.
Since $F$ was taken arbitrary, from \cite[Lemma 2.3]{DS} it follows that $Z$ is $\mathcal V$-independent.
\end{proof}

From the last claim we conclude that $|Z|=|\kappa\times(\omega_1\setminus\omega)|=\kappa$.
\end{proof}

\begin{lemma}
\label{product:of:separable:metric:groups}
Assume that 
$\kappa$
is a cardinal, $\{H_n:n\in\mathbb{N}\}$ is a family of separable metric groups in a variety $\mathcal V$ and $\{\sigma_n:n\in\mathbb{N}\}$ is a sequence of cardinals such that:
\begin{itemize}
\item[(i)] $r_\mathcal V(H_n)\ge\omega$ for every $n\in\mathbb{N}$,
\item[(ii)] $\sigma=\sup\{\sigma_n:n\in\mathbb{N}\}\ge\omega_1$,
\item[(iii)] $\Ps(\kappa,\sigma)$ holds.
\end{itemize}
Then $\prod_{n\in\mathbb{N}} H_n^{\sigma_n}$ has a $G_\delta$-dense $\mathcal V$-independent subset of size  $\kappa$. 
\end{lemma}

\begin{proof}
Define 
$$
S=\{n\in\mathbb{N}:\sigma_n\ge\omega_1\},\ \   
G=\prod_{n\in\mathbb{N}\setminus S} H_n^{\sigma_n}\ \
\mbox{ and }\ \ 
H=\prod_{n\in S} H_n^{\sigma_n}.
$$
From items (i) and (ii) of our lemma it follows that 
\begin{align}
\label{H:cong}
H\cong \prod_{i\in I} H_i', 
\mbox{ where }
|I|=\sigma
\mbox{ and }
&
\mbox{ each }
H_i'
\mbox{ is a separable metric group}
\\
&
\mbox{ satisfying }
r_\mathcal V(H'_i)\ge \omega,
\notag
\end{align}
where $\cong$ denotes the isomorphism between topological groups.  Since $|\sigma_n\times\omega_1|=\sigma_n$ for every $n\in S$, we have
$$
H^{\omega_1}\cong \prod_{n\in S}\left(H_n^{\sigma_n}\right)^{\omega_1}
\cong 
\prod_{n\in S}H_n^{\sigma_n\times\omega_1}
\cong 
\prod_{n\in S}H_n^{\sigma_n}
\cong
H.
$$
In particular, 
$$
\prod_{n\in\mathbb{N}} H_n^{\sigma_n}=G\times H\cong G\times H^{\omega_1}.
$$
Therefore, the conclusion of our lemma would follow from that of Lemma \ref{set:theoretic:lemma:1} so long as we prove that $G$ and 
$H$ satisfy the assumptions of Lemma \ref{set:theoretic:lemma:1}.
From (ii), (iii) and Proposition \ref{DS}(b) one concludes 
that $\kappa\ge\mathfrak c\ge\omega_1$.

Let us check that the assumption of item (a) of Lemma \ref{set:theoretic:lemma:1} holds. From (\ref{H:cong}) and Lemma 
\ref{big:rank:of:products} we get $r_\mathcal V(H)\ge 2^\sigma$. Since $\Ps(\kappa,\sigma)$ holds by item (iii), we have 
$2^\sigma\ge\kappa$ by Proposition \ref{DS}(c). This shows that  $r_\mathcal V(H)\ge\kappa$. 

Let us check that the assumption of item (b) of Lemma \ref{set:theoretic:lemma:1} holds. Recalling (\ref{H:cong}), we conclude that 
$$
H^\omega\cong \prod_{i\in I} \left(H_i'\right)^\omega,
\mbox{ where each }
\left(H_i'\right)^\omega
\mbox{ is a separable metric space}.
$$
Since $|I|=\sigma$ by \eqref{H:cong},
and $\Ps(\kappa, \sigma)$ holds by item (iii), Lemma \ref{Ps:in:products} allows us to conclude that $H^\omega$ has $G_\delta$-dense subset of size at most $\kappa$.

Let us check that the assumption of item (c) of Lemma \ref{set:theoretic:lemma:1} holds. Since $\sigma_n\le \omega$ for every $n\in\mathbb{N}\setminus S$, 
$G$ is a separable metric group, and so $|G|\le\cont$. Since $\Ps(\kappa,\sigma)$ holds, $\cont\le\kappa$ by Proposition \ref{DS}(b), and so $G$ itself is a $G_\delta$-dense subset of $G$ of size at most $\kappa$.
\end{proof}

\begin{corollary}
\label{Corollary:Dima}
Let $\mathbb{P}$ be the set of prime numbers and  $\{\sigma_p:p\in\mathbb{P}\}$ a sequence of cardinals such that $\sigma=\sup\{\sigma_p:p\in \mathbb{P}\}\ge\omega_1$.
If 
$\kappa$
is a cardinal such that $\Ps(\kappa,\sigma)$ holds, then  the group
\begin{equation}
\label{equation:defining:K}
K=\prod_{p\in \mathbb{P} }\mathbb{Z}_p^{\sigma_p}
\end{equation}
 contains a $G_\delta$-dense free subgroup $F$ such that $|F|=\kappa$.
\end{corollary}
\begin{proof}
Since $r(\mathbb{Z}_p)\ge\omega$ for every $p\in\mathbb{P}$, applying Lemma \ref{product:of:separable:metric:groups} 
with $\mathcal V=\mathcal A$ we can find a $G_\delta$-dense $\mathcal A$-independent subset $X$ of $K$ of size $\kappa$. 
Since $\mathcal A$-independence coincides with the usual independence for abelian groups, the subgroup $F$ of $K$ generated by $X$ is free.
Clearly, $|F|=\kappa$. Since $X\subseteq F\subseteq K$ and $X$ is $G_\delta$-dense in $K$, so is $F$.
\end{proof}

As an application, we obtain the following particular case of \cite[Lemma 4.3]{DS}.

\begin{corollary}
\label{Corollary:Dima2}
Let $\kappa$ and 
$\sigma\ge\omega_1$ 
be 
cardinals such that $\Ps(\kappa,\sigma)$ holds. Then for every compact metric non-torsion abelian group $H$ the group $H ^\sigma$ contains a $G_\delta$-dense free subgroup $F$ such that $|F|=\kappa$.
\end{corollary}

\begin{proof} 
Since $H$ is a compact non-torsion abelian group, $r(H)\geq \omega$. Applying
Lemma \ref{product:of:separable:metric:groups} with $\mathcal V=\mathcal A$, $\sigma_n=\sigma$ and $H_n=H$ for every $n\in\N$,
we can find a $G_\delta$-dense independent subset $X$ of $K=H^\sigma$ of size $\kappa$. Then the subgroup $F$ of $K$ generated by $X$ is free
and satisfies $|F|=\kappa$. Since $X\subseteq F\subseteq K$ and $X$ is $G_\delta$-dense in $K$, so is $F$.
\end{proof}

\section{Essential free subgroups of compact torsion-free abelian groups}
\label{section:6}

\begin{lemma}\label{absNEW:lemma} Let $K$ be a torsion-free abelian group and let $F$ be a free subgroup of $K$. Then there exists a free subgroup $F_0$ of $K$  
containing $F$ as a direct 
summand, such that:
\begin{itemize}
	\item[(a)] $F_0$ non-trivially meets every non-zero subgroup of $K$, and
	\item[(b)] $|F_0|=|K|$. 
\end{itemize}
\end{lemma}
\begin{proof}
Let $A:=K/F$ and let $\pi:K\to A$ be the canonical projection. Let $F_2$ be a free subgroup of $A$ with generators $\{g_i\}_{i\in I}$ such that $A/F_2$ is torsion. Since $\pi$ is surjective, for every $i\in I$ there exists $f_i\in K$, such that $\pi(f_i)=g_i$. Consider the subgroup $F_1$ of $K$ generated by  $\{f_i:i\in I\}$. As $\pi(F_1)=F_2$ is free, we conclude that $F_1\cap F=\{0\}$, so  $\pi\restriction_{F_1}:F_1\to F_2$ is an isomorphism. Let us see that the subgroup $F_0=F+F_1=F \oplus F_1$  has the required properties.
Indeed, it is free as $F_1\cap F=\{0\}$ and both $F,F_1$ are free. Moreover, $K/F_0\cong A/F_2$ is torsion and $F$ is a direct summand of $F_0$.  As $K/F_0$ is torsion,  $F_0$ non-trivially meets every non-zero subgroup of $K$, so (a) holds true.  
Since $K$ is torsion-free, (b) easily follows from (a).
\end{proof}

\begin{lemma}\label{NEW:lemma} Let $K$ be a compact torsion-free abelian group and let $F$ be a free subgroup of $K$. Then there exists a free essential subgroup $F_0$ of $K$ with $|F_0|=|K|$,  containing $F$ as a direct summand.  
\end{lemma}
\begin{proof}
Apply Lemma \ref{absNEW:lemma}.
\end{proof}

\begin{lemma}\label{Step3}
Suppose $\Min(\kappa,\sigma)$ holds, and let $\{\sigma_p:p\in\mathbb{P}\}$ be the sequence of cardinals witnessing $\Min(\kappa,\sigma)$.
Let $F$ be a free subgroup of the group $K$ as in (\ref{equation:defining:K}) with $|F|=\kappa$.
Then there exists a free  essential subgroup $F'$ of $K$ containing $F$ as a direct summand
such that  $|F'|=\kappa$.
\end{lemma}
\begin{proof} Let 
\begin{equation}\label{wtd}
\wtd(K) = \bigoplus_{p\in\mathbb{P}} \mathbb Z_{p}^{\sigma_p}\mbox{ and }F_* = F\cap \wtd(K).
\end{equation}
Then $F_*$ is a free subgroup of $\wtd(K)$, so applying Lemma \ref{absNEW:lemma} to the group 
$\wtd(K)$ and its subgroup $F_*$ we get a free subgroup $F^*$ of $\wtd(K)$ such that:
\begin{itemize}
\item[(i)] $F^*\supseteq  F_*$ and $F^*= F_*\oplus L$ for an appropriate subgroup $L$ of $F^*$; 
\item[(ii)] $F^*$ non-trivially meets every non-zero subgroup of $\wtd(K)$; 
\item[(iii)] $|F^*|= |\wtd(K)| \leq \kappa = |F|$. 
\end{itemize}
Obviously, (ii) yields that $F^*$ is essential in $\wtd(K)$.  As $\wtd(K)$ is essential in $K$ \cite{DPS}, we conclude that  $F^*$ is essential in $K$ as well. From (iii) we conclude that $F'=F+F^*$ is an essential subgroup of $K$ of size $\kappa$ containing $F$. Finally, from (\ref{wtd}) and (i) we get $F'=F + L$, and since $L\subseteq \wtd(K)$, we have
$$
F\cap L= F\cap \wtd(K) \cap L= F_*\cap L=\{0\}.
$$ 
Therefore, $F'= F\oplus L$ is free. 
\end{proof}

\begin{lemma}
\label{min:and:psc:give:necessary:embedding}
Let 
$\kappa$
and $\sigma\ge\omega_1$ be cardinals such that both  $\Min(\kappa,\sigma)$ and $\Ps(\kappa,\sigma)$ hold.
Then $F_\kappa$ admits a zero-dimensional minimal pseudocompact group topology of weight $\sigma$. 
\end{lemma}
\begin{proof}
Let $\{\sigma_p:p\in\mathbb{P}\}$ be a sequence of cardinals witnessing $\Min(\kappa,\sigma)$. In particular,
$\sigma=\sup\{\sigma_p:p\in \mathbb{P}\}$. Then the group $K$ as in (\ref{equation:defining:K}) is compact and zero-dimensional. Since $\sigma\ge\omega_1$ and $\Ps(\kappa,\sigma)$ holds, by Corollary \ref{Corollary:Dima} there exists a $G_\delta$-dense free subgroup $F$ of $K$ with $|F|=\kappa$. Since $\Min(\kappa,\sigma)$ holds, according to Lemma \ref{Step3} there exists a free essential subgroup $F'$ of $K$ containing $F$ with $|F'|=\kappa$.
Obviously $F'$ is also $G_\delta$-dense. By Theorem \ref{cr-theorem} 
$F'$ is pseudocompact. On the other hand, by the essentiality of $F'$ in $K$ and Theorem \ref{minimality-criterion}, the subgroup $F'$ of $K$ is also minimal. Being a subgroup of the zero-dimensional group $K$, the group $F'$ is zero-dimensional. Since $F'$ is dense in $K$, from \eqref{equation:defining:K} and \eqref{min-def} we have $w(F')=w(K)= \sup\{\sigma_p:p\in\mathbb{P}\}=\sigma$. Since $F'\cong  F_\kappa$, the subspace topology induced on $F'$ from $K$
will do the job.
\end{proof}

\begin{lemma}
\label{min:and:psc:give:necessary:embedding2}
Let $\kappa$ and $\sigma\ge\omega_1$ be cardinals such that $\kappa= 2^\sigma$. Then $F_\kappa$ admits a connected minimal pseudocompact group topology of weight $\sigma$. 
\end{lemma}

\begin{proof} The group $K={\widehat{\mathbb{Q}}}^\sigma$ is compact and connected. Since $\kappa=2^\sigma$, $\Ps(\kappa,\sigma)$ holds by Proposition \ref{DS}(d).
By Corollary \ref{Corollary:Dima2} there exists a $G_\delta$-dense free subgroup $F$ of $K$ with $|F|=\kappa$. 
According to Lemma \ref{NEW:lemma} there exists a free essential subgroup $F'$ of $K$ containing $F$ with $|F'|=|K|=\kappa$.
Obviously $F'$ is also $G_\delta$-dense. By Theorem \ref{cr-theorem}  $F'$ is pseudocompact. On the other hand, by the essentiality of $F'$ in $K$ and Theorem \ref{minimality-criterion}, the subgroup $F'$ of $K$ is also minimal.  Since $G_\delta$-dense subgroups of compact connected abelian groups are connected \cite[Fact 2.10(ii)]{DS}, 
we conclude that $F'$ is connected. Since $F'$ is dense in $K$, we have $w(F')=w(K)=\sigma$. Clearly, $F'\cong F_\kappa$ as $|F'|=|F|=2^\sigma=\kappa$.
Therefore, the subspace topology induced on $F'$ from $K$ will do the job.
\end{proof}
 
 \section{Proofs of the theorems from Section \ref{results:section}}
\label{proofs:section}

\begin{lemma}
\label{torsion:free:completion}
Let $G$ be a minimal torsion-free abelian group and $K$ its completion. Then:
\begin{itemize}
\item[(i)] $K$ is a compact torsion-free abelian group;
\item[(ii)] there exists a sequence of cardinals $\{\sigma_p:p\in\mathbb P\cup\{0\}\}$ such that
\begin{equation}
\label{compact:torsion:free:group}
K=\widehat{\mathbb Q}^{\sigma_0}\times\prod_{p\in\mathbb P}\mathbb Z_p^{\sigma_p}.
\end{equation}
\end{itemize}
\end{lemma}
\begin{proof}
(i) By the precompactness theorem of  Prodanov and Stoyanov (\cite[Theorem 2.7.7]{DPS}), 
$G$ is precompact, and so
$K$ is compact.  Let us show that $K$ is torsion-free. Let $x\in K\setminus\{0\}$. Assume that the cyclic group $Z=\langle x\rangle$ generated by $x$ is finite. Then $Z$ is closed in $K$ and non-trivial. Since $G$ is essential in $K$ by Theorem \ref{minimality-criterion}, it follows that $Z\cap G\not=\{0\}$. Choose $y\in Z\cap G\not=\{0\}$. Since $Z$ is finite, $y$ must be a torsion element, in contradiction with the fact that $G$ is torsion-free. 

(ii) Since $K$ is torsion-free by item (i), the Pontryagin dual of $K$ is divisible.
Now the conclusion of item (ii) of our lemma follows from \cite[Theorem 25.8]{HRoss}.
\end{proof}

\begin{proof}[\bf Proof of Theorem \ref{Step0}]
Let $K$ be the compact completion 
of $G$. 
Let $\sigma=w(K)=w(G)$. Then clearly
\begin{equation}\label{LASTeq}
|G|\leq|K|=2^\sigma.
\end{equation} 
If $\sigma=\omega$, then $|G|\leq|K|=2^\sigma=\cont$. Hence $\Min(|G|,\sigma)$ holds according to Example \ref{example:small:Stoyanov}. 
Therefore, we assume $\sigma > \omega$ for the rest of the proof.  

We consider first the case when $G$ is torsion-free. Although this part of the proof is not used in the second part covering the general case,
we prefer to include it because this provides a self-contained proof of  Theorems \ref{MPs}, \ref{continuum} and \ref{connminpsc} which concern 
only free (hence, torsion-free) groups. Let $\{\sigma_p: p\in\mathbb{P}\cup\{0\}\}$ be the sequence from the conclusion of Lemma \ref{torsion:free:completion}(ii).
Clearly,  our assumption $\sigma > \omega$ implies that $\sigma_p > \omega$ for some $p\in\mathbb{P}\cup\{0\}\}$. Hence
$\sigma=\sup\{\sigma_p: p\in\mathbb{P}\cup\{0\}\}$. 
Since $G$ is both dense and essential in $K$, from \cite[Theorems 3.12 and 3.14]{BD} we get  
$$
\sup_{p\in\{0\}\cup\mathbb P}2^{\sigma_p}\le |G|.
$$
Therefore $\Min(|G|,\sigma)$ holds in view of (\ref{LASTeq}). Since $\sigma=w(G)$, we are done.

\smallskip

In the general case, we consider the connected component $c(K)$ of $K$ and the totally disconnected quotient $K/c(K)$. Then 
$$
K/c(K)\cong \prod _{p\in\mathbb P}  K_p,
$$
 where each $K_p$ is a pro-$p$-group. Let $\sigma_p= w(K_p)$ and $\sigma_0=w(c(K))$. Our assumption $\sigma > \omega$ implies that $\sigma_p > \omega$ for some $p\in\mathbb{P}\cup\{0\}\}$, 
so that 
$$
\sigma=w(G)=w(K)=\sup_{p\in\{0\}\cup \mathbb P} \sigma_p.
$$
 By \cite[Theorems 3.12 and 3.14]{BD}, one has 
$$
|c(K)|\cdot \sup_{p\in\mathbb P} 2^{\sigma_p}\le |G|.
$$
Therefore,  
$$
\sup_{p\in\{0\}\cup \mathbb P} 2^{\sigma_p} \leq |G|\leq |K|=2^\sigma
$$
in view of (\ref{LASTeq}).
Thus $\Min(|G|,\sigma)$ holds. Since $\sigma=w(G)$, we are done.
\end{proof}

\begin{proof}[\bf Proof of Theorem \ref{weight:minimal:group}] Let $G$ be a  minimal abelian group with $w(G)\geq \kappa$.  Define $\sigma=w(G)$. Then $\Min(|G|,\sigma)$  holds by Theorem \ref{Step0}. Let $\{\sigma_n:n\in\mathbb{N}\}$ be a sequence of cardinals  witnessing $\Min(|G|,\sigma)$. That is,
\begin{equation}
\label{|G|:equation}
\sigma=\sup_{n\in\N}\sigma_n \ \ \mbox{ and }\ \ 
\sup_{n\in\N} 2^{\sigma_n}\le |G|\le 2^\sigma.
\end{equation} 
If  $\cf(\sigma)>\omega$,  then $|G|=2^{\sigma}\geq 2^{\kappa}$ by Example  \ref{exp-non-exp}(c). Assume that $\cf(\sigma)=\omega$. If $\sigma_n=\sigma$ for some $n\in\N$, then 
$|G|=2^{\sigma}\geq 2^{\kappa}$  by Example  \ref{exp-non-exp}(b). So we may additionally assume that $\sigma_n<\sigma$ for every $n\in \N$. Since $\cf(\kappa)>\omega=\cf(\sigma)$, our hypothesis  $\sigma\geq\kappa$ gives $\sigma>\kappa$.  Then $\sigma_n\ge \kappa$ for some $n\in\N$, and so $|G|\ge 2^{\sigma_n}\ge 2^{\kappa}$ by (\ref{|G|:equation}).
\end{proof}

\begin{proof}[\bf Proof of Theorem \ref{w=log||}]
By Theorem \ref{Step0}, $\Min(|G|,w(G))$ holds. Since $|G|$ is assumed to  be non-exponential, the conclusion now follows from  Proposition 
\ref{uniq}(b).
\end{proof}

\begin{proof}[\bf Proof of Theorem \ref{psc:topologies:from:minimal:ones}] Since $|F_\kappa|=\kappa$, from our assumption and Theorem \ref{Step0} we conclude that  $\Min(\kappa,\sigma)$ holds. Lemma \ref{lemma:4:8} yields that $\Ps(\kappa,\sigma)$ holds as well. Since $\sigma$ is infinite and not a strong limit, it follows that $\sigma\ge\omega_1$. Now Lemma \ref{min:and:psc:give:necessary:embedding} applies.
\end{proof}

\begin{proof}[\bf Proof of Theorem \ref{MPs}] 
The implications  (c)$\Rightarrow$(b) and (b)$\Rightarrow$(a) are  obvious.

\smallskip

(a)$\Rightarrow$(c)
Assume that $\tau_1$ is a minimal topology of weight $\sigma$ on $F_\kappa$.  Then $\sigma\ge\omega_1$ as $\kappa > \cont$. 
According to Theorem \ref{Step0} $\Min(\kappa,\sigma)$ holds.
Now assume that $\tau_2$ is a minimal topology of weight $\lambda$ on $F_\kappa$. According to Theorem \ref{CRThm} $\Ps(\kappa,\lambda)$ holds. Now Lemma \ref{Step1} yields that also $\Ps(\kappa,\sigma)$ holds true. Finally, the application of Lemma \ref{min:and:psc:give:necessary:embedding} finishes the proof.
\end{proof}

\begin{remark}
It is clear from the above proof that the topologies from items (b) and (c) of Theorem \ref{MPs} can be chosen to have the same weight $\sigma$ as the minimal topology from item (a) of this theorem.
\end{remark}

\begin{proof}[\bf Proof of Theorem \ref{continuum}] The implications (b) $\Rightarrow$(a)  and (c)$\Rightarrow$(a) are obvious.

\smallskip

(a)$\Rightarrow$(d)
Suppose that $F_\cont$ admits a minimal pseudocompact group topology. Since $F_\cont$ is free, $F_\cont$  does not
admit any compact group topology, and so $\cont=|F_\cont|\ge 2^{\omega_1}$ by Corollary \ref{min:psc:imply:compact:for:small:groups}. The converse inequality $\cont\le 2^{\omega_1}$ is clear.

\smallskip

(d)$\Rightarrow$(b) Follows from $\mathfrak c=2^{\omega_1}$ and Lemma \ref{min:and:psc:give:necessary:embedding2}. 
\smallskip

(d)$\Rightarrow$(c) Follows from $\mathfrak c=2^{\omega_1}$ and Lemma \ref{min:and:psc:give:necessary:embedding}, 
as 
$\Min(\mathfrak c,\omega_1)$ holds by Example \ref{exp-non-exp}(a),
and $\Ps(\mathfrak c,\omega_1)$ holds by Proposition \ref{DS}(a). 
\end{proof}

\begin{proof}[\bf Proof of Theorem \ref{connminpsc}]
(a)$\Rightarrow$(b) is  obvious. 

\smallskip

(b)$\Rightarrow$(c) Assume that  $\tau_1$ is a connected minimal group topology on $F_\kappa$ with $w(F_\kappa, \tau_1)=\sigma$.
Then the  completion $K$ of $(F_\kappa,\tau_1)$ satisfies the conclusion of Lemma \ref{torsion:free:completion}(ii).
Moreover, $K$ is connected. Since the zero-dimensional group 
$$
L=\prod_{p\in\mathbb P}\mathbb Z_p^{\sigma_p} 
$$
from \eqref{compact:torsion:free:group} is a continuous image of the connected group $K$, we must have $L=\{0\}$. It follows that $K=\widehat{\mathbb Q}^{\sigma_0}$.
Note that $\sigma_0=w(K)=w(F_\kappa,\tau_1)=\sigma$. That is,
$K=\widehat{\mathbb Q}^{\sigma}$.
Since $F_\kappa$ is both dense and essential in $K$ by Theorem \ref{minimality-criterion},
from \cite[Theorems 3.12 and 3.14]{BD} we get $2^{\sigma}\leq |F_\kappa|\leq |K| =2^{\sigma}$. Hence $\kappa=2^\sigma$.
\smallskip

(c)$\Rightarrow$(a) Follows from  $\kappa=2^\sigma$ and Lemma \ref{min:and:psc:give:necessary:embedding2}. 
\end{proof}

\begin{proof}[\bf Proof of Theorem \ref{no:locally:connected}] Let $G$ be a locally connected minimal abelian group and
$K$
its completion. Let $U$ be a non-empty open connected subset of $G$.  Choose an open subset $V$ of $K$ such that $V\cap G=U$. Since $U$ is dense in $V$ and $U$ is connected, so is $V$. Therefore, $K$ is locally connected. Applying Lemma \ref{torsion:free:completion}(i), we conclude that $K$ is compact and torsion-free. From \cite[Corollary 8.8]{DS} we get $K=\{0\}$. Hence $G$ is trivial as well. 
\end{proof}

\section{Final remarks and open questions}
\label{section:8}

The divisible abelian groups that admit a minimal group topology were described in \cite{D3}. Here we need only the part of this characterization for divisible abelian groups of size $\geq \mathfrak c$.

\begin{theorem}\label{divisible}\emph{\cite{D3}}
A divisible abelian group of cardinality at least $\cont$ admits some minimal group topology precisely when it admits a compact group topology. 
\end{theorem}

The concept of pseudocompactness generalizes compactness from a different  angle than that of  minimality. It is therefore quite surprising that minimality and pseudocompactness \emph{combined together} ``yield'' compactness in the class of divisible abelian groups. This should be compared with  Corollary \ref{min:psc:imply:compact:for:small:groups}, where a similar  phenomenon (i.e., minimal and pseudocompact topologizations imply compact topologization) occurs for all ``small'' groups.

The next theorem shows that the counterpart of the simultaneous minimal and pseudocompact topologization of divisible abelian groups is much easier than that of free abelian groups.

\begin{theorem}\label{div}  A divisible abelian group admits a minimal group topology and a pseudocompact group topology if and only it admits a compact group topology.
\end{theorem}

\begin{proof} The necessity is obvious.  To prove the sufficiency, suppose that a  divisible  abelian group $G$ admits both a minimal group topology and a pseudocompact group topology.  If $G$ is finite, then $G$ admits a compact group topology. If $G$ is infinite, then $|G|\geq\mathfrak c$ by Theorem \ref{van:Douwen}. Now the conclusion follows from Theorem \ref{divisible}.
\end{proof}

Our next example demonstrates that both the restriction on the cardinality in Theorem \ref{divisible} and the hypothesis of the existence of a pseudocompact group topology in Theorem \ref{div} are needed: 

\begin{example}
\label{divisible:example}
\begin{itemize}
	\item[(a)]The divisible abelian group $\mathbb Q/\mathbb Z$ admits a minimal group topology \cite{DP1}, but does not admit a pseudocompact group topology (Theorem \ref{van:Douwen}).
	\item[(b)] The divisible abelian group $\mathbb Q^{(\mathfrak c)} \oplus (\mathbb Q/\mathbb Z)^{(\omega)}$ admits a (connected) pseudocompact group topology \cite{DS}, but does not admit any minimal group topology. The latter conclusion follows from Theorem \ref{divisible} and the fact that this group does not admit any compact group topology \cite{HRoss}.
\end{itemize}
\end{example}

Let us briefly discuss the possibilities to extend our results for free abelian groups to the case of torsion-free abelian groups. Theorem \ref{div} shows that for divisible torsion-free  abelian groups the situation is in some sense similar to that of free abelian groups described in Theorem \ref{MPs}: in both cases the existence of a pseudocompact group topology and a minimal group topology is equivalent to the existence of a minimal pseudocompact (actually, compact) group topology. Nevertheless, there is a substantial difference, because free abelian groups admit no compact group topology. Another important difference between both cases is that Problem \ref{PM} is still open for torsion-free abelian groups \cite{D}: 
 
\begin{problem}
Characterize the minimal torsion-free abelian groups.
\end{problem}

A quotient of a minimal group need not be minimal even in the abelian case. This justified the isolation in \cite{DP1} of the smaller class of totally minimal groups: 

\begin{deff}
A topological group $G$ is called \emph{totally minimal} if every  Hausdorff quotient group of $G$ is minimal. Equivalently,  a Hausdorff topological group $G$ is totally minimal if  every continuous group homomorphism $f:G\to H$ of $G$ onto  a  Hausdorff  topological group $H$ is open.
\end{deff}

It is clear that compact $\Rightarrow$ totally minimal $\Rightarrow$ minimal. Therefore, Theorem \ref{MPs} makes it natural to ask the following question:

\begin{question}
Let $\kappa> \cont$ be a cardinal. 
\begin{itemize}
	\item[(a)] When does $F_\kappa$ admit a totally minimal group topology?
	\item[(b)] When does $F_\kappa$ admit a totally minimal pseudocompact group topology?
	\end{itemize}
\end{question}

More specifically, one can ask: 

\begin{question} Let $\kappa>\cont$ be a cardinal. Is the condition ``$F_\kappa$ admits a zero-dimensional totally minimal pseudocompact group topology'' equivalent to those of Theorem \ref{MPs}?
\end{question}

Since $F_\mathfrak c$ admits a totally minimal group topology \cite{P5} and a pseudocompact group topology \cite{DS}, the obvious counter-part of Theorem  \ref{continuum} suggests itself:

\begin{question}
Assume the Lusin's Hypothesis $2^{\omega_1}=\mathfrak c$.
\begin{itemize}
\item[(i)]   Does $F_\cont$ admit a totally minimal pseudocompact group topology?
\item[(ii)]  Does $F_\cont$ admit a totally minimal pseudocompact connected group topology?
\item[(iii)] Does $F_\cont$ admit a totally minimal pseudocompact zero-dimensional group topology?
\end{itemize}
\end{question}

Let us mention another class of abelian groups where both problems (Problem \ref{PM} for minimal group topologies \cite{DP2} and its counterpart for pseudocompact group topologies \cite{CRe,DS}) are completely resolved. These are the torsion abelian groups. Nevertheless, we do not know the answer of the following question: 

\begin{question} Let $G$ be a torsion abelian group that admits a  minimal group topology  and a pseudocompact group topology. Does $G$ admit also a minimal pseudocompact group topology? 
\end{question}

We finish with the question about (non-abelian) free groups. We note that the topology from Theorem \ref{minimality:of:free:groups} is even totally minimal. Furthermore, a free group $F$ admits a pseudocompact group topology if and only if $\Ps(|F|,\sigma)$ holds for some infinite cardinal $\sigma$ \cite{DS}. This justifies our final

\begin{question}
Let $F$ be a free group that admits a pseudocompact group topology.
\begin{itemize}
\item[(i)]   Does $F$ have a minimal pseudocompact group topology?
\item[(ii)]  Does $F$ have a totally minimal pseudocompact group topology?
\item[(iii)] Does $F$ have a (totally) minimal pseudocompact connected group topology?
\item[(iv)] Does $F$ have a (totally) minimal pseudocompact zero-dimensional group topology?
\end{itemize}
\end{question}

\end{document}